\documentclass[12pt]{article}
\usepackage{ytableau}
\usepackage{mathrsfs}
\usepackage{stmaryrd}
\usepackage{amsfonts,amsmath,amssymb,amscd}
\usepackage{shadow}
\usepackage{graphicx}
\usepackage{color}
\usepackage{longtable}
\usepackage{pstricks,multido}
\usepackage{hyperref}
\usepackage{qtree}
\allowdisplaybreaks

\newtheorem{thm}{Theorem}[section]
\newtheorem{lem}[thm]{Lemma}

\newtheorem{prop}[thm]{Proposition}

\newtheorem{example}{Example}[section]

\newcommand\urc{\operatorname{urc}}
\newcommand\sign{\operatorname{sign}}
\newcommand\RR{\operatorname{RR}}
\newcommand\RC{\operatorname{RC}}
\newcommand\URC{\operatorname{URC}}
\newcommand\URR{\operatorname{URR}}

\newcommand\urr{\operatorname{urr}}
\newcommand\signb{\operatorname{sign_B}}

\newcommand\nwnm{\operatorname{\overline{wm}}}
\newcommand\WM{\operatorname{WM}}
\newcommand\NWM{\operatorname{\overline{W}M}}
\newcommand\WNM{\operatorname{W\overline{M}}}
\newcommand\NWNM{\operatorname{\overline{WM}}}
\newcommand\URRB{\operatorname{URR_B}}
\newcommand\URCB{\operatorname{URC_B}}
\newcommand\urB{\operatorname{ur_B}}
\newcommand\URB{\operatorname{UR_B}}
\newcommand\topone{\operatorname{topone}}
\newcommand\wt{\operatorname{wt}}

\def\qed{\hfill \rule{4pt}{7pt}}

\def\pf{\noindent {\it{Proof.} \hskip 2pt}}

\numberwithin{equation}{section}
\numberwithin{figure}{section}

\begin{document}
\begin{center}
{\large\bf On the Sign-imbalance of Permutation Tableaux}
\end{center}

\begin{center}
Joanna N. Chen$^1$, Robin D.P. Zhou$^2$

$^{1}$College of Science\\
Tianjin University of Technology\\
Tianjin 300384, P.R. China

$^2$College of Mathematics
Physics and Information\\
Shaoxing University\\
Shaoxing 312000, P.R. China

$^1$joannachen@tjut.edu.cn,
 $^2$zhoudapao@mail.nankai.edu.cn.

\end{center}

\begin{abstract}
Permutation tableaux  were  introduced by Steingr\'{\i}msson and Williams.
 Corteel and  Kim defined the sign of a permutation tableau
 in terms of the number of
unrestricted columns. The sign-imbalance of permutation tableaux of length $n$
is the sum of  signs over  permutation tableaux of length $n$.
They  have obtained a formula for the sign-imbalance
of permutation tableaux of length $n$ by using generating functions
 and asked for a combinatorial proof. Moreover,
they raised the
question of finding a sign-imbalance formula
for type $B$ permutation tableaux  introduced by Lam and Williams.
We define a statistic $\nwnm$ over permutations and  show that  the number of unrestricted columns over permutation tableaux
of length $n$ is equally distributed with $\nwnm$ over permutations of length $n$.
  This leads to  a combinatorial interpretation of the formula of Corteel and Kim.
 For type $B$ permutation tableaux, we
 define the sign of a type $B$ permutation tableau in term of the number of
 certain rows and columns. On the other hand, we construct a bijection
 between the type $B$ permutation tableaux of length $n$ and symmetric permutations
 of length $2n$ and
 we show that the statistic $\nwnm$
 over symmetric permutations of length $2n$ is equally distributed with
 the number of certain rows and columns over type $B$ permutation tableaux of length $n$.
 Based on this correspondence and an involution on symmetric permutation of length $2n$, we obtain
a sign-imbalance formula for type $B$ permutation tableaux.

\end{abstract}



\section{Introduction}

This paper is concerned with two questions on the sign-imbalance of
permutation tableaux of type $A$ and type $B$, raised by Corteel and Kim \cite{corteel2011}.
Permutation tableaux  were introduced by Steingr\'{\i}msson and Williams \cite{Einar2006}. They are related to the enumeration of totally positive Grassmannian cells \cite{Lam2008,Postnikov2006, Scott2006, Williams2005}, as well as a statistical physics model called Partially Asymmetric Exclusion Process (PASEP) \cite{corteel,corteel2006,corteel2007,corteel2007b,corteel2011w}.
For recent studies of permutation tableaux, see, for example, \cite{alex2007,chen2011,corteel2011,corteelNadeau, Nadeau}.

A permutation tableau is defined based on the Ferrers diagram of a partition $\lambda$ for which zero parts are allowed. Let $\lambda= (\lambda_1, \lambda_2, \ldots,\lambda_r)$
be a partition, that is, $\lambda_1\geq \lambda_2\geq \cdots \geq \lambda_r \geq 0$.
The Ferrers diagram of $\lambda$ is a left-justified arrangement with $\lambda_i$
squares in the $i$th row.  The length of a Ferrers diagram is the
total number of rows and columns (including empty rows).
In particular, the length of the Ferrers diagram of the empty
partition is defined to be zero.

\begin{figure}
\begin{center}\setlength{\unitlength}{.06mm}
\begin{picture}(-800,550)(0,0)
\put(0,700){\line(1,0){500}}
\put(0,600){\line(1,0){500}}
\put(0,500){\line(1,0){500}}
\put(0,400){\line(1,0){400}}
\put(0,300){\line(1,0){200}}
\put(0,200){\line(1,0){200}}
\put(0,100){\line(1,0){100}}

\put(0,0){\line(0,1){700}}
\put(100,100){\line(0,1){600}}
\put(200,200){\line(0,1){500}}
\put(300,400){\line(0,1){300}}
\put(400,400){\line(0,1){300}}
\put(500,500){\line(0,1){200}}

\put(-80,30){\footnotesize 12}
\put(-80,130){\footnotesize10}
\put(-60,230){\footnotesize8}
\put(-60,330){\footnotesize7}
\put(-60,430){\footnotesize4}
\put(-60,530){\footnotesize2}
\put(-60,630){\footnotesize1}
\put(25,730){\footnotesize11}
\put(130,730){\footnotesize9}
\put(230,730){\footnotesize6}
\put(330,730){\footnotesize5}
\put(430,730){\footnotesize3}

\put(-1200,700){\line(1,0){500}}
\put(-1200,600){\line(1,0){500}}
\put(-1200,500){\line(1,0){500}}
\put(-1200,400){\line(1,0){400}}
\put(-1200,300){\line(1,0){200}}
\put(-1200,200){\line(1,0){200}}
\put(-1200,100){\line(1,0){100}}

\put(-1200,0){\line(0,1){700}}
\put(-1100,100){\line(0,1){600}}
\put(-1000,200){\line(0,1){500}}
\put(-900,400){\line(0,1){300}}
\put(-800,400){\line(0,1){300}}
\put(-700,500){\line(0,1){200}}

\put(-1195,10){\footnotesize 12}
\put(-1165,55){\footnotesize 11}
\put(-1095,110){\footnotesize 10}
\put(-1045,155){\footnotesize 9}
\put(-995,230){\footnotesize 8}
\put(-995,320){\footnotesize 7}
\put(-945,350){\footnotesize 6}
\put(-860,350){\footnotesize 5}
\put(-795,420){\footnotesize 4}
\put(-750,455){\footnotesize 3}
\put(-685,530){\footnotesize 2}
\put(-685,630){\footnotesize 1}

\end{picture}
\caption{The labeling of a Ferrers diagram}
\label{fig:rep}
\end{center}
\end{figure}
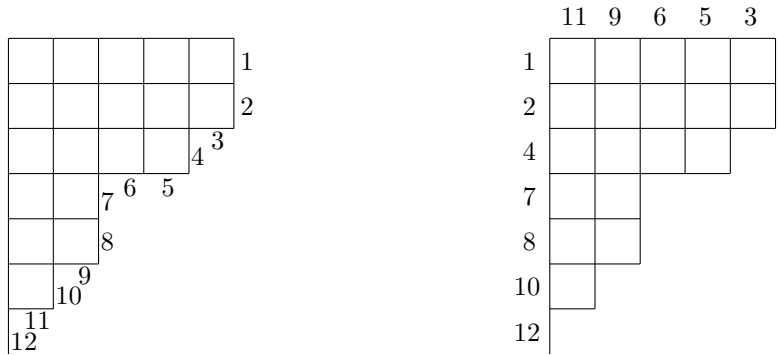

Given a Ferrers diagram $F$ of
length $n$, we label the rows and columns of $F$
as follows. First, we give labels to the steps in the south-east border with $1,2,\ldots,n$
from north-east to south-west. Then we label a row (resp.~column) with $i$ if the row
(resp.~column) contains the south (resp.~west) step with label $i$.
Notice that we may  place a row label to the left of the
first column  and place a column label  at the top of the first row,  see  Figure \ref{fig:rep}.
A row labeled with $i$ is called row $i$ and a column labeled with $j$ is called
column $j$.
We use $(i, j)$ to denote the cell in row $i$ and column $j$.

For a partition $\lambda$, a permutation tableau of shape $\lambda$ is a $0,1$-filling of the Ferrers diagram of $\lambda$  satisfying the following conditions:
\begin{enumerate}
\item Each column has at least one $1$;
\item There is no $0$ with a $1$ above (in the same column) and a $1$ to
  the left (in the same row).
\end{enumerate}

The length of a permutation tableau is defined to be the length of the corresponding
Ferrers diagram.  Denote by $\mathcal{PT}(n)$ the set of permutation tableaux of length $n$.
Figure \ref{fig:permalt} illustrates a permutation tableau
of length $12$.

\begin{figure}
\begin{center}\setlength{\unitlength}{.06mm}
\begin{picture}(550,550)(-115,0)
\put(0,700){\line(1,0){500}}
\put(0,600){\line(1,0){500}}
\put(0,500){\line(1,0){500}}
\put(0,400){\line(1,0){400}}
\put(0,300){\line(1,0){200}}
\put(0,200){\line(1,0){200}}
\put(0,100){\line(1,0){100}}

\put(0,0){\line(0,1){700}}
\put(100,100){\line(0,1){600}}
\put(200,200){\line(0,1){500}}
\put(300,400){\line(0,1){300}}
\put(400,400){\line(0,1){300}}
\put(500,500){\line(0,1){200}}

\put(-80,30){$12$}
\put(-80,130){$10$}
\put(-60,230){$8$}
\put(-60,330){$7$}
\put(-60,430){$4$}
\put(-60,530){$2$}
\put(-60,630){$1$}
\put(15,730){$11$}
\put(120,730){$9$}
\put(220,730){$6$}
\put(320,730){$5$}
\put(420,730){$3$}

\put(30,630){$0$}
\put(30,530){$0$}
\put(30,430){$1$}
\put(30,330){$0$}
\put(30,230){$0$}
\put(30,130){$1$}

\put(130,630){$1$}
\put(130,530){$0$}
\put(130,430){$1$}
\put(130,330){$1$}
\put(130,230){$1$}

\put(230,630){$0$}
\put(230,530){$1$}
\put(230,430){$1$}

\put(330,630){$0$}
\put(330,530){$0$}
\put(330,430){$1$}

\put(430,630){$1$}
\put(430,530){$1$}

\end{picture}
\caption{A permutation tableau $T$ of length $12$ }
\label{fig:permalt}
\end{center}
\end{figure}

In their study of combinatorics of permutation tableaux in connection with PASEP,
Corteel and  Williams \cite{corteel2007} introduced the concepts of a row-restricted $0$ and an unrestricted row.
A $0$ in a permutation tableau is said to be row-restricted if there is a $1$ above
(in the same column). A row is called unrestricted if it does not contain any row-restricted
$0$. Otherwise, it is called a restricted row.
Let $T$ be a permutation tableau  with $k$ columns, and let $\urr(T)$ be the number of unrestricted rows of $T$.
Besides, Corteel and  Williams \cite{corteel2007}
introduced the weight $\wt(T)$ as the total number of $1$'s in $T$ minus $k$,
and used the notation $\topone(T)$ for the number of $1$'s in the first row of $T$. They defined the polynomial
$$
F_{\lambda, \alpha, \beta}(q)=
\sum_{T} q^{\wt(T)} \alpha^{-\topone(T)}\beta^{-\urr(T)+1},
$$
where the sum ranges over  permutation tableaux $T$ of shape $\lambda$.
Using the matrix ansatz for the PASEP model, they derived a formula for $F_{\lambda, \alpha, \beta}(q)$.
Corteel and Nadeau \cite{corteelNadeau}
obtained an explicit formula for the generating function of permutation tableaux  of length $n$ with respect to the statistics $\urr$ and $\topone$:
\begin{equation}
\label{urrtop}
  \sum_{T\in\mathcal{PT}(n)} x^{\urr(T)-1} y^{\topone(T)} = (x+y)_{n-1},
\end{equation}
where $(x)_n$ denote the rising factorial, that is, $(x)_0=1$ and
  $(x)_n=x(x+1)\cdots(x+n-1)$ for $n\geq1$.

Corteel and Kim \cite{corteel2011} gave two bijective proofs of (\ref{urrtop}).
Furthermore, they introduced the concepts of a column-restricted $0$ and an unrestricted column.
A $0$ in a permutation tableau is called column-restricted if there is a $1$ to the left
(in the same row). In the same vain, one can define unrestricted columns and restricted columns. For a permutation tableau $T$,  let $\urc(T)$ denote the number of unrestricted columns of $T$. Corteel and Kim obtained the following generating function
of permutation tableaux  of  length $n$ with respect to the statistic $\urc$.

\begin{thm}\label{th:gf}
We have
$$\sum_{n\geq0} \sum_{T\in \mathcal{PT}(n)} t^{\urc(T)} x^n = \frac{1+E_t(x)}{1+(t-1)x E_t(x)},$$
where
$$E_t(x)=\sum_{n\geq1} n (t)_{n-1} x^n.$$
\end{thm}

The above generating function leads to a formula for
the sign-imbalance  of permutation tableaux of length $n$.
For a permutation tableau $T$,  the sign of $T$ is given by
$\sign(T)=(-1)^{\urc(T)}$. Let
\begin{equation}\label{sn}
 s(n)=\sum_{T\in\mathcal{PT}(n)} \sign(T).
\end{equation}
 Setting $t=-1$ in Theorem \ref{th:gf}, Corteel and Kim
derived the following formula for   $s(n)$.

\begin{thm}
\label{th:c1}
Assume that  $n=4k+r$, where $0 \leq r <4$. Then
\begin{equation*}
s(n) = \left\{
  \begin{array}{ll}
(-1)^k \cdot 2^{2k}, & \mbox{if $r=0$ or $r=1$,}\\[6pt]
0, & \mbox{if $r=2$,}\\[6pt]
(-1)^{k+1} \cdot 2^{2k+1}, & \mbox{if $r=3$.}
 \end{array} \right.
\end{equation*}
\end{thm}

Corteel and Kim \cite{corteel2011}
asked for a combinatorial proof
of Theorem \ref{th:c1}.
In answer to this question, we introduce a permutation
statistic  $\nwnm$.
More precisely, for a permutation $\pi$,
we define $\nwnm(\pi)$  to be the number of
element $\pi_i$ such that $\pi_i < i$ and $\pi_i$
does not appear in the middle of a decreasing subsequence of length three.
We show that the statistic $\urc$  over permutation tableaux
of length $n$ is equally distributed with the statistic $\nwnm$ over permutations of length $n$. Moreover, for $n\geq 4$, we build
a parity reversing involution on the set of permutations $\pi=\pi_1 \pi_2 \cdots \pi_n$ such that $\pi_1 \pi_2 \pi_3 \pi_4 \neq 1342, 1432, 2341 , 2431$.
 Using this involution, we obtain a combinatorial proof of Theorem \ref{th:c1}.

Moreover, the construction of the aforementioned involution
implies the following recurrence relation for $s(n)$.

\begin{thm} \label{recurrence}
For $n \geq 3$,
\[s(n)=2 s(n-1)- 2s(n-2).\]
\end{thm}

Notice that the above recurrence relation along with the
initial values $s(1)=1$ and $s(2)=0$ also leads to  a
proof of Theorem \ref{th:c1}.

The second result of this paper is concerned with a question proposed
by Corteel and Kim on permutation tableaux of type $B$.
Type $B$ permutation tableaux were
introduced by Lam and Williams \cite{Lam2008}, and
further studied by  Corteel and Kim \cite{corteel2011},
and  Corteel, Josuat-verg\`{e}s and Kim \cite{corteel2012}.
A type $B$ permutation tableau is defined based on a shifted Ferrers diagram.
Let $F$ be a Ferrers diagram with $k$ columns.
Unlike the underlying Ferrers diagram of a permutation tableau,
for the type $B$ case, both empty rows and empty columns are allowed
in a  Ferrers diagram.
The shifted Ferrers diagram
of $F$, denoted by $\bar{F}$, is defined to be the diagram obtained from $F$
by adding $k$ rows of size $1,2,\ldots,k$ at the top of the diagram.
\begin{figure}
\begin{center}\setlength{\unitlength}{.053mm}
\begin{picture}(-800,700)(0,0)
\put(0,900){\line(1,0){100}}
\put(0,800){\line(1,0){200}}
\put(0,700){\line(1,0){300}}
\put(0,600){\line(1,0){400}}
\put(0,500){\line(1,0){500}}
\put(0,400){\line(1,0){500}}
\put(0,300){\line(1,0){400}}
\put(0,200){\line(1,0){200}}
\put(0,100){\line(1,0){200}}
\put(0,100){\line(1,0){100}}
\put(0,0){\line(1,0){100}}

\put(0,0){\line(0,1){900}}
\put(100,0){\line(0,1){900}}
\put(200,100){\line(0,1){700}}
\put(300,300){\line(0,1){400}}
\put(400,300){\line(0,1){300}}
\put(500,400){\line(0,1){100}}

\put(-110,830){$-9$}
\put(-110,730){$-7$}
\put(-110,630){$-4$}
\put(-110,530){$-3$}
\put(-110,430){$-1$}
\put(-64,330){$2$}
\put(-64,230){$5$}
\put(-64,130){$6$}
\put(-64,30){$8$}

\put(15,930){$9$}
\put(120,930){$7$}
\put(220,930){$4$}
\put(320,930){$3$}
\put(420,930){$1$}

\put(25,830){$\ast$}
\put(125,730){$\ast$}
\put(225,630){$\ast$}
\put(325,530){$\ast$}
\put(425,430){$\ast$}

\put(-1200,400){\line(1,0){500}}
\put(-1200,300){\line(1,0){400}}
\put(-1200,200){\line(1,0){200}}
\put(-1200,100){\line(1,0){200}}
\put(-1200,100){\line(1,0){100}}
\put(-1200,0){\line(1,0){100}}

\put(-1200,0){\line(0,1){400}}
\put(-1100,0){\line(0,1){400}}
\put(-1000,100){\line(0,1){300}}
\put(-900,300){\line(0,1){100}}
\put(-800,300){\line(0,1){100}}

\end{picture}
\caption{A Ferrers diagram and its labeled corresponding shifted Ferrers diagram}
\label{fig:shifted}
\end{center}
\end{figure}
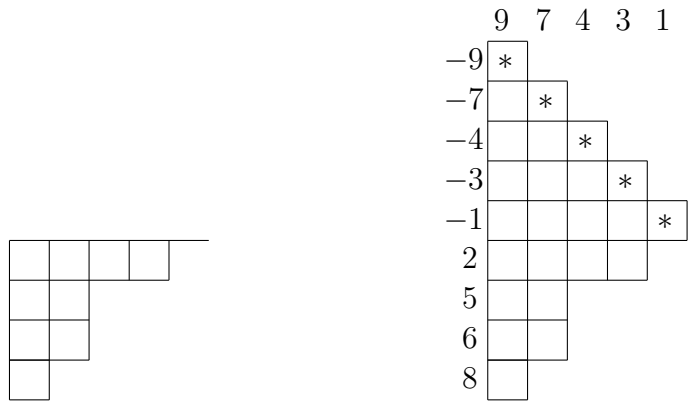

For a Ferrers diagram $F$, the length of the corresponding shifted Ferrers diagram $\bar{F}$
is defined to be the length of $F$,  and the diagonal of $\bar{F}$ is the set of rightmost cells of the added rows.
A diagonal cell is a cell on the diagonal.
We label the added row by the opposite number of the label of the column in which
the rightmost cell in the row locates.
The labels of the other rows and columns remain the same with $F$. Figure \ref{fig:shifted} illustrates  a  Ferrers diagram and its corresponding labeled shifted Ferrers diagram, where the diagonal cells are marked with stars.

A type $B$ permutation tableau is a $0,1$-filling of a
shifted Ferrers diagram satisfying the following conditions:
\begin{enumerate}
  \item each column has at least one $1$.
  \item there is no $0$ which has a $1$ above (in the same column) and a $1$ to the
        left (in the same row).
  \item if a $0$ is in a diagonal cell, then it does not have any $1$ to the left
        (in the same row).
\end{enumerate}
The length of a type $B$ permutation tableau is defined to be the length of the corresponding
shifted Ferrers diagrams.
Let $\mathcal{PT}_B(n)$ be the set of type $B$ permutation tableaux of length $n$.
Figure \ref{2.1fig} gives  a type $B$ permutation tableau of length 8.

\begin{figure}[h]
\begin{center}
\ytableausetup
{mathmode, boxsize=1.5em}
 \begin{ytableau}
     \none&\none[8]&\none[6]&\none[3]&\none[2]\\
     \none[-8]& 0\\
     \none[-6]&1& 1 \\
     \none[-3]&0& 0& 0\\
     \none[-2]&0& 1& 0& 1\\
     \none[1]&1& 1& 1 &1\\
     \none[4]&0& 1\\
     \none[5]&0& 0\\
     \none[7]& 1
 \end{ytableau}
 \end{center}
\caption{A type $B$ permutation tableau of length $8$.}
\label{2.1fig}
\end{figure}

Corteel and  Kim \cite{corteel2011} raised the question of finding a sign-imbalance formula for type $B$ permutation tableaux. In answer to this question, we introduce the concepts of unrestricted rows and
columns for a type $B$ permutation tableau.
Then we define
the sign of a type $B$ permutation tableau by the total
number of unrestricted rows and columns.
We show that there is a bijection between type $B$ permutation tableaux
of length $n$ and symmetric permutations on $\{1, 2, \ldots, 2n\}$,
where a symmetric
permutation $\pi$ of length $2n$
is a permutation on $[2n]$ such that $\pi_i+\pi_{2n+1-i}=2n+1$
for $1 \leq i \leq n$.
Using this correspondence, we derive a sign-imbalance formula for type $B$ permutation tableaux of length $n$.

\section{A Combinatorial Proof of Theorem \ref{th:c1} }\label{sec:typeA}

In this section, we introduce a permutation statistic
$\nwnm$. We show that the statistic $\nwnm$ over
permutations of length $n$ is equally distributed
with the statistic $\urc$ over $\mathcal{PT}(n)$. For $n \geq 4$, we exhibit a parity reversing involution on the set of permutations with the first four elements not equal to $1342$, $1432$, $2341$ nor $2431$. Based on this involution, we give a combinatorial proof of Theorem \ref{th:c1}.
Moreover, we derive a recurrence relation for $s(n)$ as given
in Theorem \ref{recurrence}.
By taking the initial values of $s(n)$ into consideration,
we give another proof of Theorem \ref{th:c1}.

We first introduce the  permutation statistic $\nwnm$.
 Let $[n]=\{1,2, \ldots, n \}$ and
 $S_n$ be the set of permutations on $[n]$.
 Given a permutation $\pi=\pi_1 \pi_2 \cdots \pi_n$ of length $n$, the index $i$ is said to be a  weak excedance  of  $\pi$
if $\pi_i \geq i$. Otherwise, it is called a non-weak excedance.
A permutation $\pi$ is said to contain a pattern $\tau$ if there exists a subsequence of $\pi$ that has the same relative order as $\tau$.
Otherwise, $\pi$ is said to avoid  $\tau$.
The element $\pi_i$ is called a
mid-point of $\pi$ if it is the middle point of a decreasing subsequence of
length three of $\pi$,
namely, there exist $j<i$ and $k>i$ such that $\pi_j>\pi_i >\pi_k$.
Otherwise, $\pi_i$ is called a non-mid-point.

For a permutation $\pi=\pi_1 \pi_2 \cdots \pi_n$, we define the following four sets.
\begin{eqnarray*}
\WM(\pi) &=&\{\pi_i| ~i\ \text{is a weak excedance and a mid-point of $\pi$} \},\\[5pt]
\NWM(\pi)&=& \{\pi_i|~ i\ \text{is a non-weak excedance and a  mid-point of $\pi$ }   \},\\[5pt]
\WNM(\pi)&=& \{\pi_i| ~i\ \text {is a weak excedance and a non-mid-point of $\pi$ }  \},\\[5pt]
\NWNM(\pi)&=& \{\pi_i|~ i \ \text {is a non-weak excedance and a non-mid-point of $\pi$ }  \}.
\end{eqnarray*}
Set $\nwnm(\pi)$ to be the number of elements in $\NWNM(\pi)$.
As an example, given a permutation $\pi=6, 5, 1, 10, 4, 3, 8, 9, 2, 11, 7, 12$ of length $12$,
we have $ \NWNM(\pi)=\{1,2,7 \}$ and $\nwnm=3$.

In the following, we aim to show that the statistic $\nwnm$ over $S_n$ is equally distributed with the statistic $\urc$ over $\mathcal{PT}(n)$.
To achieve this, we  first recall a bijection $\Phi$ from $\mathcal{PT}(n)$ to $S_n$, which was given by Steingr\'{\i}msson and Williams \cite{Einar2006}.

 A zigzag path on a permutation tableau $T \in \mathcal{PT}(n)$ is a path entering from the left of a row or the top of a column, going to the east or to the south changing the direction alternatively whenever it meets a $1$ until exiting the tableau. For convenience, we denote the
 path entering from $i$ by $P_i$.
  Then $\Phi$ is defined to be the permutation $\pi=\pi_1 \pi_2 \cdots \pi_n$ where $\pi_i=j$ if the zigzag path $P_i$ exits $T$ from row $j$ or column $j$. As an example, for the permutation tableau $T$ given in Figure \ref{fig:permalt},
   $\Phi(T)=6,5,1,10,4,3,8,9,2,11,7,12$.

The map $\Phi$ has the following
properties, of which the proof of the former one  can be found in Steingr\'{\i}msson and Williams \cite{Einar2006},
and hence, we omit it here.

\begin{prop}
\label{row-wex}
Given $T \in \mathcal{PT}(n)$, let $\pi=\Phi(T)$.
If $i$ is a weak excedance of $\pi$, then $i$ is precisely a row label of $T$. Otherwise,
$i$ is a column label. In particular, $\pi$ is a permutation
with $k$ weak excedance if and only if $T$ has exactly $k$ rows.
\end{prop}

\begin{prop}
\label{intersect}
Let $T$ be a permutation tableau in $\mathcal{PT}(n)$. For $1 \leq i <j \leq n$, if the zigzag paths $P_i$ and $P_j$ intersect, they can only intersect at points, not edges. Moreover, the intersecting points, except the first one,  must correspond to a
$1$ in $T$.

\end{prop}
\pf
Without loss of generality, we may assume that $i$ is a row label and $j$ is a column label of $T$. The other cases can be proved similarly. If paths $P_i$ and $P_j$ intersect,
let $x$ be the first point they meet. It's easily seen that
path $P_i$ travels to $x$ from west, while path
$P_j$ travels to $x$ from north.
To give a proof of this proposition,
we consider two cases.

If $x$ corresponds to a $1$ in $T$, when meeting $x$,
path $P_i$ turns to south and path $P_j$ turns to
east. Hence, we see that paths $P_i$ and $P_j$ intersect at the
point $x$, not an edge starting from $x$. Let $y$ be the second point they meet, if there exists. Clearly, path $P_i$ travels to $y$ from west, while path
$P_j$ travels to $y$ from north.  We claim that $y$ corresponds to a $1$ in $T$. Otherwise, there will exist a $0$ corresponding to $y$ which has a $1$ above it and a $1$ to its left, which contradicts to the definition of the permutation
tableau.  Then, by a similar analysis, we may prove that  paths $P_i$ and $P_j$ intersect at the point $y$, not an edge starting from $y$.
Using induction on the number of the
intersecting points, we see that in this case paths $P_i$ and
$P_j$ can only intersect at points and each intersecting point corresponds to a $1$ in $T$.

  If $x$ corresponds to a $0$ in $T$, when meeting $x$, path $P_i$ goes through $x$ to east, while path $P_j$ goes through $x$ to south. Hence paths $P_i$ and $P_j$  intersect at  the point $x$, not an edge starting from $x$. Set $y$ to be the second point they meet, if there exists.
  It is routine to check that $P_i$ travels to $y$ from north, while $P_j$ travels to $y$
  from west. By a similar analysis with the above case, we see that $P_i$ and $P_j$ intersect
  at the point $y$, which corresponds to a $1$ in $T$. Using induction on the number of the intersecting points, we see that  paths $P_i$ and
  $P_j$ can only intersect at points and each intersecting point, except the first point, corresponds to a $1$ in $T$.
  Combining the above two cases, we complete the proof. \qed

 Now we are ready to give the proof of the equidistribution of the statistics $\urc$ and $\nwnm$.
Given a permutation tableau $T$, let $\RR(T),\URR(T),$ $\RC(T), \URC(T)$  denote the set
of the restricted rows, unrestricted rows, restricted columns and unrestricted columns of $T$, respectively. Clearly, we have $\urc(T)=|\URC(T)|$ and $\urr(T)=|\URR(T)|$.
As an example, for the permutation tableau $T$ given in Figure \ref{fig:permalt},
we have $\URC(T)=\{3,9,11\}$ and $\urc(T)=3$.

We will show that $\nwnm$ is equally distributed with $\urc$ by exhibiting a
one-to-one correspondence between $\RC(T)$ and $\NWM(\pi)$ for a permutation tableau $T$ and $\pi=\Phi(T)$,
as given in the following lemma.

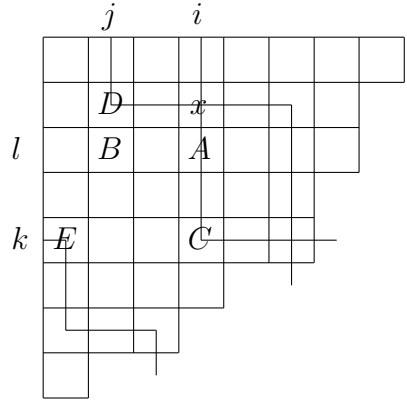
\begin{figure}
\begin{center}\setlength{\unitlength}{.06mm}
\begin{picture}(550,550)(-40,0)
\put(0,800){\line(1,0){800}}
\put(0,700){\line(1,0){800}}
\put(0,600){\line(1,0){700}}
\put(0,500){\line(1,0){700}}
\put(0,400){\line(1,0){600}}
\put(0,300){\line(1,0){600}}
\put(0,200){\line(1,0){400}}
\put(0,100){\line(1,0){300}}
\put(0,0){\line(1,0){100}}

\put(0,0){\line(0,1){800}}
\put(100,0){\line(0,1){800}}
\put(200,100){\line(0,1){700}}
\put(300,100){\line(0,1){700}}
\put(400,200){\line(0,1){600}}
\put(500,300){\line(0,1){500}}
\put(600,300){\line(0,1){500}}
\put(700,500){\line(0,1){300}}
\put(800,700){\line(0,1){100}}

\put(20,330){$E$}
\put(-70,330){$k$}
\put(120,630){$D$}
\put(130,830){$j$}
\put(120,530){$B$}
\put(-70,530){$l$}
\put(320,530){$A$}
\put(320,330){$C$}
\put(330,830){$i$}
\put(325,630){$x$}

\put(0,350){\line(1,0){50}}
\put(50,350){\line(0,-1){200}}
\put(50,150){\line(1,0){200}}
\put(250,150){\line(0,-1){100}}

\put(150,650){\line(0,1){150}}
\put(150,650){\line(1,0){400}}
\put(550,650){\line(0,-1){300}}
\put(550,350){\line(1,0){100}}

\put(350,350){\line(0,1){450}}
\put(350,350){\line(1,0){200}}
\put(550,350){\line(0,-1){100}}
\end{picture}
\caption{The description of $P_i$, $P_j$ and $P_k$ }
\label{fig:zigzag}
\end{center}
\end{figure}

\begin{lem}\label{lem:eqdistri}
Given a permutation tableau $T \in \mathcal{PT}(n)$, let $\pi=\pi_1\pi_2\ldots \pi_n=\Phi(T)$.
 There is a bijection between
  $\URC(T)$  (resp.  $\RC(T)$, $\RR(T)$, $\URR(T))$
and $\NWNM(\pi)$
(resp. $\NWM(\pi)$, $\WM(\pi)$, $\WNM(\pi))$.
\end{lem}
\pf
By Proposition \ref{row-wex}, we deduce that
$|\URC(T)$$ \cup \RC(T)|=|\NWNM(\pi)\cup \NWM(\pi)|.$
Hence, to prove that there is a bijection between $\URC(T)$ and $\NWNM(\pi)$, it suffices
to show that there is a bijection between $\RC(T)$ and $\NWM(\pi)$.
Similarly, to prove that there is a bijection between $\URR(T)$ and $\WNM(\pi)$,
it suffices to prove that there is a bijection  between $\RR(T)$ and $\WM(\pi)$.
Here, we will provide the bijection   between $\RC(T)$ and $\NWM(\pi)$ only, as the bijection  between $\RR(T)$ and $\WM(\pi)$ can be given similarly.

We claim that $\pi$ is a
bijection from $\RC(T)$ to $\NWM(\pi)$, namely, $i\in \RC(T)$ if and only if $\pi_i\in \NWM(\pi)$.

On the one hand, we proceed to show that for any $i\in \RC(T)$, we have $\pi_i\in \NWM(\pi)$.
Given $i\in \RC(T)$,
it follows from  Proposition \ref{row-wex} that $i$ is a non-weak excedance.
Hence, to prove $\pi_i \in \NWM(\pi)$, it is left to prove that $\pi_i$ is a mid-point. Since column $i$ is restricted, there is at least one column-restricted $0$ in it. Let $A$ be the lowermost column-restricted $0$ in column $i$, which is located in row $l$, see Figure \ref{fig:zigzag} as a description.  Since $A$ is column-restricted, there must be a $1$
to the left of $A$. Denote the nearest $1$ to the left of $A$ by $B$ and suppose $B$ is in column $j$, where $j>i$. It is easy to see that
in column $i$ there is no $1$ above $A$.
By the definition of permutation tableau, each column contains
at least one $1$.  So let $C$ be the
topmost $1$ in column $i$ and  assume that
$C$ is in row $k$, where $k <i$.
Clearly, we have $k< i<j$, to prove that $\pi_i$ is a mid-point, we aim to prove that
$\pi_k > \pi_i >\pi_j$.

 Firstly, we prove that $\pi_j < \pi_i$.
We claim that the first intersecting point of the zigzag paths $P_i$ and $P_j$ corresponds to a $0$ in $T$. Let $D$ be the topmost $1$ in column $j$, which is not lower than $B$. Path $P_j$ travels to $D$  and then turn east. Path $P_i$ travels to $C$  and then turn east. Since $B$ is the nearest $1$ to the left of $A$,
the element in cells $(l,g)$, where $i< g <j$, are all column-restricted $0$'s.
Hence, all the cells above
row $l$ and between column $j$ and column $i$ are filled by $0$'s.
This means that  the first intersecting point of $P_i$ and $P_j$, which is denoted by $x$, corresponds to a $0$ in column $i$. The claim is verified.
Moreover, by Proposition \ref{intersect}, the following
intersecting points, if there exist, must correspond to $1$'s
in $T$.
Then it is not hard to see that after the intersecting
point $x$, path $P_j$ is always on the upper right of
$P_i$. It follows that  $\pi_j <\pi_i$, as desired.

Next, we show that $\pi_i<\pi_k$. Assume that $E$ is the
leftmost $1$ in row $k$. Clearly, it is not to the right of $C$. If $P_k$ and $P_j$ do not intersect, then $P_i$ is always on the upper right of $P_k$. It is obvious that $\pi_i<\pi_k$.
Otherwise, suppose the first intersecting point of $P_k$ and $P_j$ is $y$. We claim that $y$ corresponds to a $1$ in $T$.
If point $E$ coincides with $C$, then clearly $P_k$ and $P_j$
intersect at $C$, which corresponds to a $1$, as desired.
If point $E$ is to the left of $C$, $P_k$ travels to $E$ and then turn south. Then it is not hard to see that  path $P_i$  travels to $y$ from north, and path $P_k$ travels to $y$ from west. Hence, $y$
can not be a $0$. Otherwise, $y$ will have a $1$ above it and a $1$ to its left, a contradiction.  The claim is verified. Again, by Proposition \ref{intersect},  the following
intersecting points, if there exist, must correspond to $1$'s
in $T$.  This means that $P_i$ is always on the upper right of $P_k$. It follows that $\pi_i<\pi_k$, as desired.

On the other hand,  we need to show that for each $\pi_i\in \NWM(\pi)$,  we have $i\in \RC(T)$.
Since $\pi_i$ is a mid-point of $\pi$, then there exist $\pi_k$ and $\pi_j$ such that $k<i<j$ and $\pi_k>\pi_i>\pi_j$.
In view of $\pi_i < i$, we deduce that $\pi_j <j$.
Let $T=\Phi^{-1}(\pi)$, by Proposition \ref{row-wex}, $i$ and $j$ are both column labels of $T$. If there is no column-restricted zero in column $i$, then the topmost $1$ of column $i$ can't be lower than the topmost $1$ of column $j$. By a similar discussion as above, we can obtain $\pi_i<\pi_j$, a contradiction. It follows that
$i \in \RC(T)$. This completes the proof. \qed

\begin{example}
Let $T$ be the permutation tableau given by Figure  \ref{fig:permalt}. Then,
$\pi=\Phi(T)=6, 5, 1, 10, 4, 3, 8, 9, 2, 11, 7, 12$. We have
\begin{align*}
\centering
\RR(T)=\{2,7,8\}&\longleftrightarrow  \WM(\pi)=\{5,  8, 9 \},\\[5pt]
\RC(T)=\{5,6\} &\longleftrightarrow  \NWM(\pi)=\{ 3, 4 \},\\[5pt]
\URR(T)=\{1,4,10,12\}&\longleftrightarrow  \WNM(\pi)=\{6, 10, 11, 12\},\\[5pt]
\URC(T)=\{3,9,11\}&\longleftrightarrow  \NWNM(\pi)=\{1, 2, 7\}.
\end{align*}
\end{example}

As a consequence of Lemma \ref{lem:eqdistri},
we see that the statistic $\urc$ is equidistributed with the statistic $\nwnm$.
It follows that
\begin{equation*}
s(n)=\sum_{T\in \mathcal{PT}(n) }(-1)^{\urc(T)}=\sum_{\pi\in S_n }(-1)^{\nwnm(\pi)}.
\end{equation*}

To give a proof of Theorem
\ref{th:c1}, we proceed to exhibit a parity reversing involution on a subset of $S_n$,
which are stated in the following lemma.

\begin{lem}\label{allinv}
There is an involution  $\varphi$ on the subset $W_n$ of $S_n$, where
\[W_n=\{\pi \in S_n | ~\pi_1 \pi_2 \pi_3 \pi_4 \neq 1342, 1432, 2341, 2431\}.\]
Moreover, $\nwnm(\pi)$ and $\nwnm(\varphi(\pi))$ have different parties for $\pi \in W_n$.
\end{lem}

Before presenting the proof of Lemma \ref{allinv}, we  first construct an involution on a subset of $W_n$.

\begin{lem}
\label{lem:12}
There is an involution $\phi$ on the subset $V_n$ of $W_n$, where
\[V_n=\{\pi \in S_n | ~\pi_1\neq 1,2\} \cup  \{\pi \in S_n | ~\pi_1 \pi_2=12  \text{ or } \pi_1\pi_2=21\}.\]
Moreover,  $\nwnm(\pi)$ and $\nwnm(\phi(\pi))$ have different parties for $\pi \in V_n$.
\end{lem}

\pf Given a permutation $\pi$ in $V_n$, let $\tau=\phi(\pi)$
be given by exchanging the positions of $1$ and $2$ in $\pi$.
It is easy to check that $\phi$ is an involution on $V_n$. In the following, we proceed to show that $\phi$ is parity reversing, namely, $\nwnm(\pi)$ and $\nwnm(\tau)$ have different parities.

We may assume that $1$ is to the left of $2$ in $\pi$. For $1 \leq i \leq n$, if $\pi_i \neq 1$ and $\pi_i \neq 2$, it is easy to check that $\pi_i \in \NWNM (\pi)$ if and only if $\tau_i \in \NWNM(\tau)$.
Hence, to prove this lemma, we need only to consider the elements $1$ and $2$ for the two cases below.

\begin{itemize}
  \item The case $\pi_1\neq 1,2$.

  Since $\pi_1\neq 1,2$ and $1$ is to the left of $2$ in $\pi$, we see that $\pi_2 \neq 2$.
 This implies that $2$ is a non-weak excedance of $\pi$. By the fact that
 $1$ is to the left of $2$, we see that $2$ is a non-mid point of $\pi$. It follows that
$2 \in \NWNM(\pi)$.
Since $2$ is to the left of $1$ in $\tau$,
 then $\tau_1 21$ forms a $321$-pattern of $\tau$. It follows that $2$ is a mid-point of $\tau$, which implies that  $2 \notin \NWNM(\tau)$. It is not hard to check  that $1 \in \NWNM(\pi)$ and $1 \in \NWNM(\tau)$.
 We conclude  that $\nwnm(\tau)=\nwnm(\pi)-1.$  This completes the proof of this case.

  \item The case $\pi_1\pi_2=12$.

  It is easily checked that
$1,2 \notin \NWNM(\pi)$, $2 \notin \NWNM(\tau)$, while $1 \in \NWNM(\tau)$.
This means that $\nwnm(\tau)=\nwnm(\pi)+1.$ Hence, $\nwnm(\pi)$ and $\nwnm(\tau)$ have different parities.
\end{itemize}

Combing the above two cases, we complete the proof.\qed

Now, we proceed to give the proof of Lemma \ref{allinv}.

\noindent
{\it Proof of Lemma \ref{allinv}.}
By Lemma \ref{lem:12}, we have given an involution $\phi$ on $V_n$. Hence, to prove this lemma, it is left to construct a parity reversing involution $\theta$
on $U_n$, where $U_n=W_n/V_n$. Clearly, we have $U_n=\{\pi \in W_n |~\pi_1=1 \text{ and } \pi_2 \neq 2\} \cup \{\pi \in W_n |~\pi_1=2 \text{ and } \pi_2 \neq 1\}$.
Without loss of generality,  we assume that $\pi_1=1$. The proof of the
cases that $\pi_1=2$ can be performed similarly.

Given a permutation $\pi$ in $U_n$ with $\pi_1=1$,
we proceed to construct its image under the map $\theta$.
If  $\pi_1 \pi_2 \pi_3 \pi_4=1324$, let $\tau=\theta(\pi)$, where $\tau$ is obtained from $\pi$ by exchanging
 the positions of $3$ and $4$. It is easily checked that in this case we have $\theta^2(\pi)=\pi$ and
$\nwnm(\pi)=\nwnm(\tau)-1$.

If  $\pi_1 \pi_2 \pi_3 \pi_4 \neq 1324$ and $\pi_1 \pi_2 \pi_3 \pi_4 \neq 1423$, let
$\pi_i \pi_j \pi_k$ be the subsequence of $\pi$ containing $2,3$ and $4$.
Set $\tau=\theta(\pi)$,
where $\tau$ is obtained from $\pi$ by exchanging the positions of
$\pi_j$ and $\pi_k$. Clearly, in this case we have  $\theta^2(\pi)=\pi$. In the following, we aim to compute the relations between $\nwnm(\pi)$ and
 $\nwnm(\tau)$. To achieve this, we consider two cases.
\begin{itemize}
  \item The case $j>3$.

  We may assume that $\pi_j < \pi_k$. By the definition of $\theta$, we see that $\tau_i \tau_j \tau_k= \pi_i \pi_k \pi_j$. By a routine analysis, we
  can obtain that $\pi_j \in \NWNM(\pi)$, $\pi_k \in \NWNM(\pi)$ and $\pi_j \in \NWNM(\tau)$. If $i > 2$, then $\tau_2 \pi_k \pi_j$ forms a $321$ pattern of $\tau$.  If $i = 2$, then $\tau_3 \pi_k \pi_j$ forms a $321$ pattern of $\tau$. Thus, we conclude that $\pi_k \notin \NWNM(\tau)$.
  It follows that in this case we have
  $\nwnm(\pi)=\nwnm(\tau)+1$.

  \item The case $j=3$.

  By the definition of $U_n$ and the assumption that
  $\pi_1\pi_2\pi_3 \pi_4 \neq 1324,1423$, we deduce that $i=2$ and $k>4$.
  We may assume that $\pi_j < \pi_k$. Since $\pi_2 \neq 2$, we have
  $\pi_j=2$. Then it is routine to check that
  $\pi_j \in \NWNM(\pi)$, $\pi_k \in \NWNM(\pi)$ and $\pi_j \in \NWNM(\tau)$. By the fact that
  $\pi_k \geq 3$, we see that $\tau_3=\pi_k$ is a weak excedance of $\tau$. Hence,   $\pi_k \notin \nwnm(\tau)$. It follows that in this case we have
  $\nwnm(\pi)=\nwnm(\tau)+1$.

\end{itemize}
Combing the above two cases, we complete the proof. \qed

Based on Lemma \ref{lem:eqdistri} and Lemma \ref{allinv}, we are ready to give a proof of Theorem \ref{th:c1}.

\noindent
{\it Proof of Theorem \ref{th:c1}.}
Set $n=4k+r(0 \leq r <4)$.
 Following from Lemma \ref{lem:eqdistri}, it suffices  to prove that
\begin{equation*}
s(n)=\sum_{\pi\in S_n }(-1)^{\nwnm(\pi)}= \left\{
  \begin{array}{ll}
(-1)^k \cdot 2^{2k}, & \mbox{if $r=0$ or $r=1$,}\\[6pt]
0, & \mbox{if $r=2$,}\\[6pt]
(-1)^{k+1} \cdot 2^{2k+1}, & \mbox{if $r=3$.}
 \end{array} \right.
\end{equation*}

Let $R_n$ be the set of permutation $\pi=\pi_1\pi_2 \cdots \pi_n$ such that  $\pi_{4i-3}\pi_{4i-2}\pi_{4i-1}\pi_{4i}$ is a permutation of  $\{4i-3,4i-2,4i-1,4i\}$  that is order isomorphic to $1342,1432,2341$ or $2431$  for $1 \leq i \leq k$.
In the following, we proceed to construct a parity reversing involution $\chi$ on
$S_n/R_n$.

Assume that $\pi$ is a permutation in $S_n/R_n$. By the definition of $R_n$, we see that there exist integers $i(1 \leq i \leq k)$ such that $\pi_{4i-3}\pi_{4i-2}\pi_{4i-1}\pi_{4i}$ is not a permutation of  $\{4i-3,4i-2,4i-1,4i\}$ that is order isomorphic to $1342,1432,2341$ or $2431$.
Let $j$ be the minimum element among these integers.
Suppose that $p$ is the permutation obtained from $\pi$ by deleting the first
$4j-4$ elements and subtracting each remaining element by $4j-4$. Clearly, $p$
is a permutation in $W_m$, where $m=4(k-j)+4+r$. Write $q=\theta(p)$.
We define $\chi$ by letting $\tau=\tau_1 \tau_2 \cdots \tau_{n}=\chi(\pi)$, where $\tau_1\tau_2 \cdots \tau_{4j-4}=\pi_1 \pi_2 \cdots \pi_{4j-4}$ and $\tau_i=q_{i-4j+4}+4j-4 ~(i \geq 4j-4)$. By the fact that $\theta$ is a parity reversing involution, it is routine to check that $\chi$ is a parity reversing involution.
It follows that \[\sum_{\pi\in S_n/R_n }(-1)^{\nwnm(\pi)}=0.\]
Thus, we have $s(n)=\sum_{\pi\in R_n }(-1)^{\nwnm(\pi)}$.
To compute $s(n)$, we consider the following four cases.
\begin{itemize}
\item  The case $r=0$, namely, $n=4k$

Assume that $\pi$ is a permutation in $ R_{4k}$.
For $1 \leq i \leq k$, we have $\{\pi_i,\pi_{i+1},\pi_{i+2},\pi_{i+3}\}$
is order isomorphic to $1342,1432,2341$ or $2431$.
Then, it is not hard to check that there exists exactly one element in $\{\pi_i,\pi_{i+1},\pi_{i+2},\pi_{i+3}\}$ which is contained in the set $\NWNM(\pi)$. It follows that
 $\nwnm(\pi)=(-1)^k$. Since $|R_{4k}|=2^{2k}$, we have $s(n)=(-1)^k \cdot 2^{2k}$ for $n=4k$.

\item  The case $r=1$, namely, $n=4k+1$

Given a permutation $\pi$ in $R_{4k+1}$,
by the definition of $R_{4k+1}$, we see that $\pi_n=n$.
Clearly, $n \notin \NWNM(\pi)$.
We deduce that $s(n)=(-1)^k \cdot 2^{2k}$  for $n=4k+1$.

\item  The case $r=2$, namely, $n=4k+2$

Given a  permutation $\pi$ in $R_{4k+2}$,
by the definition of $R_{4k+2}$, we see that
$\pi_{n-1}\pi_n$ is a permutation of the set $\{n-1,n\}$.
Assume that $\pi_{n-1}=n-1$ and $\pi_{n}=n$. Setting $\tau$ to be the permutation obtained from $\pi$ by exchanging the positions of $n-1$ and $n$. It is easily seen that $\nwnm(\pi)=\nwnm(\tau)-1$.
It follows that $s(n)=0$ for $n=4k+2$.

\item  The case $r=3$, namely, $n=4k+3$

Let $\pi$ be the permutation in $R_{4k+3}$. By the definition of $R_{4k+3}$, we see that $\pi_{n-2}\pi_{n-1}\pi_n$ is a permutation of $\{n-2,n-1,n\}$.
If $\pi_{n-2}\pi_{n-1}\pi_n$ is order isomorphic to $123$ or $321$, let $\tau$ be the permutation obtained from $\pi$ by exchanging the positions of $n-2$ and $n-1$. It is routine to check that $\nwnm(\pi)=\nwnm(\tau)-1$.
Thus, s(n) equals to the sign-imbalance of the set
$R_{4k+3}'$, which is a  subset of $R_{4k+3}$ with $\pi_{n-2}\pi_{n-1}\pi_n$ order isomorphic to $132$ or $231$.  Notice that $\nwnm(\pi)=(-1)^{k+1}$ for each $\pi \in R_{4k+3}'$ and
$|R_{4k+3}'|=2^{2k+1}$. We conclude that $s(n)=(-1)^{k+1} \cdot 2^{2k+1}$ for $n=4k+3$.

\end{itemize}

Combining the above four cases, we complete the proof. \qed

Up to now, we have given a combinatorial proof of Theorem
 \ref{th:c1}.  In fact, there is another partially combinatorial proof of this theorem.
We first derive the following recurrence relation for $s(n)$
as given in Theorem \ref{recurrence}.

\noindent
{\it Proof of Theorem \ref{recurrence}.}
By Lemma \ref{lem:12}, we deduce that
$\sum_{\pi \in V_n} (-1)^{\nwnm(\pi)}=0.$ Hence, to obtain a recurrence relation
for $s(n)$, it is left  to deal with the remaining cases for $\pi \in S_n/V_n$.
Clearly,   $S_n/V_n=\{\pi \in S_n |~\pi_1=1 \text{ and } \pi_2 \neq 2\} \cup \{\pi \in S_n |~\pi_1=2 \text{ and } \pi_2 \neq 1\}$.

Given a permutation $\pi$ with $\pi_1=1$ and $ \pi_2 \neq 2$,
let $\tau=\tau_1 \tau_2\cdots \tau_{n-1}$ be the permutation given by $\tau_i=\pi_{i+1}-1$.
Since   $1$ is a weak excedance of $\pi$, we see that
$1 \notin \NWNM(\pi)$. It follows that $\nwnm(\tau)=\nwnm(\pi).$
Notice that $\tau_1 \neq 1$,
 we deduce that
\begin{align*}
    \sum_{{\pi \in S_n} \atop {\pi_1=1  \text{ and }   \pi_2 \neq 2} }(-1)^{\nwnm(\pi)}&=
    \sum_{\tau \in S_{n-1},  \tau_1 \neq 1 }(-1)^{\nwnm(\tau)}\\[3pt]
    &=s(n-1)-s(n-2).
\end{align*}

Moreover, it is routine to check that
the sign-imbalance of the set $\{\pi \in S_n |~\pi_1=1 \text{ and } \pi_2 \neq 2\}$ is equal to that of the set $\{\pi \in S_n |~\pi_1=2 \text{ and } \pi_2 \neq 1\}$. Hence, we conclude that $s(n)=2s(n-1)-2s(n-2)$. This completes the proof.
\qed

Notice that $s(1)=1$ and $s(2)=0$. By Theorem \ref{recurrence}, we give another proof of Theorem  \ref{th:c1}.

\section{The sign-imbalance of permutation tableaux of type B}\label{sec:typeB}

In this section, we define the sign of a type $B$ permutation tableau $T$, which is
denoted by $\signb(T)$.
We show that there is a bijection between type $B$ permutation tableaux
of length $n$ and symmetric permutations on $[2n]$.
Using this correspondence, we derive a sign-imbalance formula for type $B$ permutation tableaux, which is given as follows.
\begin{thm}
\label{thm:sign-imbalanceB}
If $n=2k+r(0 \leq r <2)$, then
\begin{align*}
s_B(n)=\sum_{T\in\mathcal{PT}_B(n)} \signb(T)&= \left\{
  \begin{array}{ll}
2^{\,\frac{n}{2}\!}, & \mbox{if $r=0$,}\\
0, & \mbox{if $r=1$.}
 \end{array} \right.
\end{align*}
\end{thm}

We first introduce two symmetric constructions corresponding to a type $B$ permutation tableau $T$ of length $n$.
Add the cells  obtained by reflecting the non-diagonal cells in $T$
about the diagonal line to  $T$,
we get a symmetric construction $T_s$, which we call
symmetric tableau.
It is easily seen that the length of the symmetric tableau $T_s$ is $2n$.
Label the steps in the south-east border of $T_s$ with $1,2,\ldots,2n$ from
north-east to south-west. We label a row (resp.~column) of $T_s$ with $i$ if the row
(resp.~column) contains the south (resp.~west) step with label $i$.
By condition $3$ in the definition of type $B$ permutation
tableaux, we see that there is no $0$ in $T_s$ which has a $1$ above (in the same column) and a $1$ to the left (in the same row).
Hence, if we remove all the rows and columns of $T_s$ that have no $1$'s and keep the labels  unchanged,
we obtain a symmetric permutation tableau of type $A$, which is denoted $T_A$.
Here, we remark that in this paper we always use $T_s$ and $T_A$ to present the corresponding symmetric tableau and symmetric permutation tableau of $T$, respectively.
As an example, for the type $B$ permutation tableau $T$ in Figure \ref{2.1fig},
its corresponding $T_s$ and  $T_A$ are given in Figure \ref{2.2fig}.

\begin{figure}[h]
\begin{center}
\ytableausetup
{mathmode, boxsize=1.4em}
 \begin{ytableau}
     \none& \none[16]& \none[14]& \none[11]&\none[10]&\none[8]&\none[5]&\none[4]&\none[2]\\
     \none[1]&  0& 1& 0& 0& 1& 0& 0& 1\\
     \none[3]&  1& 1& 0& 1& 1& 1& 0   \\
     \none[6]&  0& 0& 0& 0& 1         \\
     \none[7]&  0& 1& 0& 1& 1         \\
     \none[9]&  1& 1& 1& 1            \\
     \none[12]& 0& 1                  \\
     \none[13]& 0& 0                  \\
     \none[15]& 1
 \end{ytableau}
 \qquad \qquad \qquad
\ytableausetup
{mathmode, boxsize=1.4em}
 \begin{ytableau}
     \none& \none[16]& \none[14]& \none[11]&\none[10]&\none[8]&\none[5]&\none[2]\\
     \none[1]& 0& 1& 0& 0& 1& 0& 1\\
     \none[3]& 1& 1& 0& 1& 1& 1  \\
     \none[6]& 0& 0& 0& 0& 1      \\
     \none[7]& 0& 1& 0& 1& 1      \\
     \none[9]& 1& 1& 1& 1         \\
     \none[12]& 0& 1\\
     \none[15]& 1
 \end{ytableau}
\end{center}
\caption{The corresponding symmetric tableau $T_s$ (left) and symmetric permutation tableau $T_A$ (right) of $T$ given in Figure \ref{2.1fig}.}
\label{2.2fig}
\end{figure}

Given a type $B$ permutation tableau $T$, we say a $0$ in $T$ is row-restricted if it has a $1$ above (in the same column).
While, a $0$ is said to be column-restricted if it
has a $1$ to the left (in the same row).
Note that in this paper a $0$ in the diagonal cell is not treated as row-restricted, which is
different from the definitions given by Corteel and Kim  \cite{corteel2012}.
We define the set of type B unrestricted rows  $\URRB(T)$
and the set of type $B$ unrestricted columns  $\URCB(T)$ of $T$ as follows.
\begin{align*}
\URRB(T)=&\{i~|~ i~ \text{is a positive row label satisfies that row $i$ contains no row-restricted $0$}\\
         &\text{and contains at least one $1$}\},\\
\URCB(T)=&\{j~|~j ~\text{is a column label satisfies that column $j$ contains no column-restricted $0$  } \\
&\text{and row $-j$ contains no row-restricted $0$ }\}.
\end{align*}
Let
$\URB(T)=\URRB(T)~ \cup ~ \URCB(T)$ and
$\urB(T)=|\URB(T)|$.
Define the sign of the  type $B$ permutation tableau $T$ by
\[\signb(T)=(-1)^{\urB(T)}.\]
 By the construction of $T_A$ from $T$, it is easy to check that $\urB(T)=\urc(T_A)$. It follows that
\[\signb(T)=(-1)^{\urc(T_A)}.\]

Let $SP_{2n}$ be the set
of the symmetric permutations on $[2n]$.
In the following, we shall construct a bijection $\Phi_B$ between
 $\mathcal{PT}_B(n)$ and  $SP_{2n}$, which allows us to translate the statistic $ur_B$ on
  $\mathcal{PT}_B(n)$ to the statistic $\nwnm$ on $SP_{2n}$.

Let $T$ be a type $B$ permutation tableau of length $n$.
We define $\pi=\Phi_B(T)$ to be the permutation on $[2n]$ which is obtained
from $T_s$ by using the zigzag map.
Let $S$ be the set of labels of the rows and columns in
$T_s$ which contain no $1$.
In fact,  $\pi=\Phi_B(T)$ can be also obtained by
computing $\pi'=\Phi(T_A)$ on $[2n]\setminus S$ and then setting
the element in $S$ to be  fixed points of $\pi$.
As an example, for the type $B$ permutation tableau $T$ given in Figure
\ref{2.1fig},
we have that $\pi=\Phi_B(T)=8,1,12,4,3,7,11,2,15,6,10,14,13,5,16,9$.
While $S=\{4,13\}$ and
\begin{equation*}
 \pi'=\left(
 \begin{array}{cccccccccccccc}
 1& 2 &3 & 5& 6& 7 &8&9&10&11&12 &14& 15& 16\\ 8&1&12&3&7&11&2&15&6&10&14&5&16&9\\
 \end{array}
\right).
\end{equation*}
It should be remarked that $\pi'(i)$ is said to be a  weak excedance if $\pi'(i) \geq i$. Otherwise, $\pi'(i)$ is said to be a non-weak excedance.

To show that $\Phi_B$ is a bijection from  $\mathcal{PT}_B(n)$ to $SP_{2n}$, we need the
following proposition.
\begin{prop}\label{pro:3.1}
Let $T \in \mathcal{PT}_B(n)$ and $\pi=\Phi_B(T)$.
If $i$ is a fixed point of $\pi$, then $i$ is a column label of
$T_s$ when $1\leq i\leq n$ and $i$ is a row label of
$T_s$ when $n+1\leq i\leq 2n$.
If $\pi_i$ is an excedance, then $i$ is a row label of $T_s$, and if $\pi_i$ is a non-weak excedance, then $i$ is a column label of $T_s$.
\end{prop}

\pf Since $T \in \mathcal{PT}_B(n)$, we see that
 there is at least one $1$ in each column of $T$.
It follows that each row in $T_s$ with label $1\leq i\leq n$
and each column in $T_s$ with label $n+1\leq j\leq 2n$ contain at least one $1$.
Hence if $i$ is a fixed point of $\pi=\Phi_B(T)$, then  $i$ is a column label of $T_s$ when $1\leq i\leq n$ and $i$ is a row label of
$T_s$ when $n+1\leq i\leq 2n$.
Since the non-fixed points of $\pi$ is the permutation $\Phi(T_A)$,
then the remaining parts of the proposition follows directly from Proposition \ref{row-wex}. The proof is completed.
\qed

Based on Proposition \ref{pro:3.1}, we have the following theorem.

\begin{thm}
$\Phi_B$ is a bijection from  $\mathcal{PT}_B(n)$ to $SP_{2n}$.
\end{thm}
\pf Given a type $B$ permutation tableau $T\in \mathcal{PT}_B(n)$, let
$\pi=\Phi_B(T)$.
We claim that $\pi$ is a symmetric permutation on $[2n]$.
Since $T_s$ is a symmetric tableau, the zigzag paths $P_i$
and $P_{2n+1-i}$ is symmetric about the diagonal line of $T_s$.
Then, it is easy to check  that $\pi_i+\pi_{2n+1-i}=2n+1$ for $1\leq i\leq 2n$.
Hence, $\pi$ is a symmetric permutation on $[2n]$, as claimed.

To prove that $\Phi_B$ is a bijection, we will give an explicit description
of its inverse.
Given  $\pi \in SP_{2n}$, set $T$ to be the type B
permutation tableau obtained by the following procedure.

First, we construct the symmetric tableau $T_s$ corresponding to $T$.
Compute the fixed points, excedances and non-weak excedances of $\pi$.
Using Proposition \ref{pro:3.1}, it is easy  to obtain a Ferrers diagram $F$, which is
the shape of $T_s$.  Hence, to get $T_s$,
we need to fill the cells of $F$ with $0$'s or $1$'s.
If $i$ is a fixed point of $\pi$,  fill the cells in row $i$ or column $i$ with $0$'s.
Denote $F_A$ the Ferrers diagram obtained by removing all the rows and columns which have been already filled. It should be noted  that
the labels of the columns and rows in $F_A$ remain the same with
that of $F$.
Let $\pi'$ be the permutation obtained from $\pi$ by deleting all the fixed points.
Let $T_A=\Phi^{-1}(\pi')$. It is not hard to check that
$F_A$ is the shape of $T_A$. Then, fill $F_A$ with $0$'s and $1$'s such that it equals to
$T_A$. Thus, we have constructed $T_s$. And the type $B$ permutation tableau
$T$ can be obtained by removing all the cells on the upper right of the diagonal line of $T_s$.
It is clear from the above construction that $T$ is a type $B$ permutation tableau, and
moreover, it is the inverse image of $\Phi_B^{-1}(\pi)$.
This completes the proof.
\qed

Now, we are ready to
translate the statistic $\urB$ to the statistic $\nwnm$.
\begin{lem}\label{lem:B}
The statistic $\urB$ on $\mathcal{PT}_B(n)$ is equaidistributed with the
statistic $\nwnm$ on $SP_{2n}$.
\end{lem}
\pf
Given $T\in \mathcal{PT}_B(n) $, we write $\pi=\Phi_B(T)$.
Recall that $\urB(T)=|\URC(T_A)|$.
Hence, to prove this lemma, we need only to show that  there is a bijection between $\URC(T_A)$ and $\NWNM(\pi)$.
By the definition of the map $\Phi_B$, we see that $\pi$
can be obtained from $\pi'=\Phi(T_A)$ by setting all the other elements in $[2n]$ as fixed points.
Following from Lemma \ref{lem:eqdistri}, there is a bijection between
$\URC(T_A)$ and $\NWNM(\pi')$. Thus, to prove this lemma, it suffices to
show that $\NWNM(\pi')=\NWNM(\pi)$.

Firstly, we claim that if $\pi_i \in \NWNM(\pi')$, then $\pi_i \in \NWNM(\pi)$.
Since $\pi_i$ is a non-weak excedance of $\pi'$, there exists $j<i$ such that
$\pi_j>\pi_i$. By the fact that $\pi_i$ is a non-mid point of $\pi'$, we deduce that
 $\pi_l >\pi_i$ for any non-fixed point $l>i$ of $\pi$. For  fixed point $k >i$ of $\pi$,  it is easy to see that $\pi_k > \pi_i$. Hence,  all the elements that are after $\pi_i$ in $\pi$ are larger than $\pi_i$.
 It follows that $\pi_i$ is a non-mid point of $\pi$.
 Since $\pi_i$ is a non-weak excedance of $\pi'$, then $\pi_i$
 is a non-weak excedance of $\pi$. Thus,
 $\pi_i \in \NWNM(\pi)$. The claim is verified.

 Notice that the fixed points of $\pi$ are weak excedance of $\pi$.
 It is easy to verify that if $\pi_i \in \NWNM(\pi)$, then
 $\pi_i \in \NWNM(\pi')$. Hence,
we conclude that $\NWNM(\pi')=\NWNM(\pi)$.
This completes the proof.
\qed

Based on Lemma \ref{lem:B}, we proceed to give a proof of Theorem \ref{thm:sign-imbalanceB}. Notice that the proof of this theorem is similar to that of Theorem \ref{th:c1}.
We just outline the main idea of the proof, details are omitted.

\noindent
{\it Proof of Theorem \ref{thm:sign-imbalanceB}.}
By Lemma \ref{lem:B}, we see that
Theorem \ref{thm:sign-imbalanceB} is equivalent to
\begin{equation}\label{sign-imbalanceB}
s_B(n)=\sum_{\pi\in SP_{2n}}(-1)^{\nwnm(\pi)}=\left\{
                                       \begin{array}{ll}
                                       2^{\frac{n}{2}}, & \hbox{if $n$ is even,}
                                         \\[5pt]
                                          0, & \hbox{if $n$ is odd.}

                                       \end{array}
                                     \right.
\end{equation}
Define a subset of $SP_{2n}$ as follows
$$
SS_{2n}=\{\pi \in SP_{2n} | \pi_1 \pi_2 \neq 12 ~\text{and} ~\pi_1
\pi_2 \neq 21 \}.
$$
For the case that $\pi \in SS_{2n}$,
exchanging the positions of $1, 2$ and the positions of  $2n-1,2n$ to get $\tau$.
It can be shown that the parities of $\nwnm(\pi)$ and $\nwnm(\tau)$ are different.
The proof of this fact can be performed similarly to that of Lemma \ref{lem:12}.
It follows that
\[\sum_{\pi \in SS_{2n}} (-1)^{\nwnm(\pi)}=0.\]
Hence, we deduce that
\begin{align*}
s_B(n)&=\sum_{\pi\in SP_{2n}/SS_{2n}}(-1)^{\nwnm(\pi)}\\[5pt]
    &=\sum_{{\pi \in SP_{2n}} \atop {\pi_1=1  \text{ and }   \pi_2= 2} }(-1)^{\nwnm(\pi)}+\sum_{{\pi \in SP_{2n}} \atop {\pi_1=2  \text{ and }   \pi_2= 1} }(-1)^{\nwnm(\pi)}\\[5pt]
    &=\sum_{\pi'\in SP_{2n-4} }(-1)^{\nwnm(\pi')}+ \sum_{\pi''\in SP_{2n-4} }(-1)^{\nwnm(\pi'')}\\[5pt]
    &=2s_B(n-2).
\end{align*}
%
%
%

Taking the initial values into consideration,
 we arrive at \eqref{sign-imbalanceB}.
This completes the proof.
\qed


\end{document}